\documentclass[11pt]{article}

\usepackage{amsthm} 
\usepackage{fullwidth}
\usepackage[top=1.5in, bottom=1.5in, left=1in, right=1in]{geometry}
\usepackage{xypic}
\usepackage{amsgen}
\usepackage{amsmath}
\usepackage{amstext}
\usepackage{amsbsy}
\usepackage{amsopn}
\usepackage{amsfonts}
\usepackage{amssymb}
\usepackage{eepic}
\usepackage[pdftex]{graphicx}
\usepackage{epsf}
\usepackage[usenames,dvipsnames]{xcolor}
\usepackage{pstricks}
\usepackage{makeidx}
\xyoption{all}

\def\mapright#1{\smash{\mathop{\longrightarrow}\limits^{#1}}}
\def\tra#1{\smash{\mathop{\mid\kern
-1pt\joinrel\relbar\joinrel\relbar}\limits^{*}_{#1}}}
\def\longtra#1{\smash{\mathop{\mid\kern
-1pt\joinrel\relbar\joinrel\relbar\joinrel\relbar}\limits^{*}_{#1}}}
\def\vlongtra#1{\smash{\mathop{\mid\kern-1pt\joinrel\relbar\joinrel\relbar\joinrel\relbar\joinrel\relbar}\limits^{*}_{#1}}}
\def\vvlongtra#1{\smash{\mathop{\mid\kern
-1pt\joinrel\relbar\joinrel\relbar\joinrel\relbar\joinrel\relbar\joinrel\relbar}\limits^{*}_{#1}}}
\def\vvvlongtra#1{\smash{\mathop{\mid\kern
-1pt\joinrel\relbar\joinrel\relbar\joinrel\relbar\joinrel\relbar\joinrel\relbar\joinrel\relbar}\limits^{*}_{#1}}}
\def\etra#1{\smash{\mathop{\mid\kern
-1pt\joinrel\relbar\joinrel\relbar}\limits_{#1}}}
\def\mapleft#1{\smash{\mathop{\longleftarrow}\limits^{#1}}}
\def\vlongrightarrow{\relbar\joinrel\longrightarrow}
\def\vvlongrightarrow{\relbar\joinrel\vlongrightarrow}
\def\vvvlongrightarrow{\relbar\joinrel\vvlongrightarrow}
\def\vvvvlongrightarrow{\relbar\joinrel\vvvlongrightarrow}
\def\vvvvvlongrightarrow{\relbar\joinrel\vvvvlongrightarrow}

\def\longmapright#1{\smash{\mathop{\vlongrightarrow}\limits^{#1}}}

\def\vvlongmapright#1{\smash{\mathop{\vvvlongrightarrow}\limits^{#1}}}

\def\vvvvlongmapright#1{\smash{\mathop{\vvvvvlongrightarrow}\limits^{#1}}}

\def\A{{\cal{A}}}

\def\iff{\Leftrightarrow}
\def\Rw{\Rightarrow}
\def\oo{\overline}

\def\Gr{\text{Gr}}
\def\B{{\cal{B}}}
\def\C{{\cal{C}}}
\def\CR{{\cal{CR}}}
\def\CS{{\cal{CS}}}

\def\I{{\cal{I}}}
\def\L{{\cal{L}}} 
\def\M{{\cal{M}}}

\def\S{{\cal{S}}}
\def\X{{\cal{X}}}

\def\fix{\mbox{Fix}}
\def\UA{{\bf UA}}
\def\UE{{\bf UE}}
\def\tak{\mbox{Tak}}

\def\gcd{\mbox{gcd}}

\def\rk{\mbox{rk}}

\def\per{\mbox{Per}}
\def\ker{\mbox{Ker}}

\def\endo{\mbox{End}}
\def\aut{\mbox{Aut}}

\def\min{\mbox{min}\,}

\def\R{\mathrel{{\mathcal R}}}

\def\N{\mathbb{N}}
\def\H{\mathrel{{\mathcal H}}}

\def\G{{\cal{G}}}
\def\V{{\cal{V}}}

\def\ZZ{\mathbb{Z}} 
\def\Q{\mathbb{Q}}

\def\p{\varphi}

\def\inv{^{-1}}

\def\J{\mathrel{{\mathcal J}}} 

\newcommand{\dotcup}{\ensuremath{\mathaccent\cdot\cup}}

\renewcommand{\geq}{\geqslant}
\renewcommand{\leq}{\leqslant}

\def\bi{\begin{itemize}}
\def\ei{\end{itemize}}
\def\ben{\begin{enumerate}}
\def\een{\end{enumerate}}
\def\beq{\begin{equation}}
\def\eeq{\end{equation}}

\newtheorem{T}{Theorem}[section]
\newcommand{\bt}{\begin{T}}
\newcommand{\et}{\end{T}}
\newcommand{\ftd}{$\square$\end{T}}

\newtheorem{Proposition}[T]{Proposition}
\newcommand{\bp}{\begin{Proposition}}
\newcommand{\ep}{\end{Proposition}}
\newcommand{\fpd}{$\square$\end{Proposition}}

\newtheorem{Definition}[T]{Definition}
\newcommand{\bd}{\begin{Definition}}
\newcommand{\ed}{\end{Definition}}

\newtheorem{Lemma}[T]{Lemma}
\newcommand{\bl}{\begin{Lemma}}
\newcommand{\el}{\end{Lemma}}
\newcommand{\fld}{$\square$\end{Lemma}}

\newtheorem{Corol}[T]{Corollary}
\newcommand{\bc}{\begin{Corol}}
\newcommand{\ec}{\end{Corol}}
\newcommand{\fcd}{$\square$\end{Corol}}

\newtheorem{Result}[T]{Result}
\newcommand{\br}{\begin{Result}}
\newcommand{\er}{\end{Result}}
\newcommand{\frd}{$\square$\end{Result}}

\newtheorem{Remark}[T]{Remark}
\newcommand{\brem}{\begin{Remark}}
\newcommand{\erem}{\end{Remark}}
\newcommand{\fremd}{$\square$\end{Remark}}

\newtheorem{Example}[T]{Example}
\newcommand{\be}{\begin{Example}}
\newcommand{\ee}{\end{Example}}

\newtheorem{Problem}[T]{Problem}
\newcommand{\bq}{\begin{Problem}}
\newcommand{\eq}{\end{Problem}}

\newcommand{\bproof}{\begin{proof}[\emph{\textbf{Proof}}]}
\newcommand{\eproof}{\end{proof}}

\newlength{\lengtha} 
\newlength{\lengthb} 

\def\abstract#1{\par\bigskip
\begingroup\small
\baselineskip=12truept
\begin{center}ABSTRACT\end{center}
\par\medskip\par\noindent
\null\hfill\hbox{\vbox{\hsize=5truein\noindent#1}}
\hfill\null\par\endgroup\par}

\title{Takahasi semigroups}
\author{{\bf M\'ario J.\ J.\ Branco}\\ 
$ $\\ {\em Departamento de Matem\'atica and CAUL/CEMAT,} \\
  {\em Faculdade de Ci\^encias, Universidade de Lisboa}\\ 
  {\em Campo Grande, 1749-016 Lisboa, Portugal}\\
  {\em email:} mjbranco@fc.ul.pt\\
$ $\\
{\bf Gracinda M.\ S.\ Gomes}\\ 
$ $\\ {\em Departamento de Matem\'atica and CAUL/CEMAT,} \\
  {\em Faculdade de Ci\^encias, Universidade de Lisboa}\\ 
  {\em Campo Grande, 1749-016 Lisboa, Portugal}\\
{\em email:} gmcunha@fc.ul.pt\\
$ $\\
{\bf Pedro V.\ Silva}\\ $ $\\
{\em Centro de
Matem\'{a}tica, Faculdade de Ci\^{e}ncias, Universidade do
Porto,}\\ {\em R.\ Campo Alegre 687, 4169-007 Porto, Portugal}\\
{\em email:} pvsilva@fc.up.pt}

\date{\today}

\begin{document}

\maketitle

\begin{center}\small

\bigskip

\end{center}

\abstract{Takahasi's theorem on chains of subgroups of bounded rank in
  a free group is generalized to several classes of semigroups. As an
  application, it is proved that the subsemigroups of 
  periodic points are finitely generated and periodic orbits are
  bounded for arbitrary endomorphisms for various semigroups.
  Some of these results feature classes such as completely
  simple semigroups, Clifford semigroups or monoids defined by
  balanced one-relator presentations. 
  In addition to the background on semigroups, 
  proofs involve arguments over groups and finite automata.
}

\section{Introduction} \label{Sec:intro}

In a paper of 1950, Takahasi proved the following, often called 
{\em Takahasi's Theorem} (the same result was proved independently by Higman 
in the paper~\cite{Hig51} of 1951): 
  
\bt
\label{taka}
{\rm \cite{Tak}}
Let $F$ be a free group and let
$K_1 \leq K_2 \leq \dotsm$ be an ascending chain of finitely generated
subgroups of $F$. If the rank of the subgroups in the chain is
bounded, then the chain is stationary.
\et

Recall that in view of Nielsen's Theorem \cite[Proposition~I.2.6]{LS},
every finitely generated subgroup of a free group is free, so in Takahasi's
Theorem {\em rank} means the cardinality of a basis. 

The concept of rank admits a natural generalization to arbitrary
groups. Given a group $G$, we define its group rank, 
denoted by $\rk_{\G}(G)$, to be the 
minimum cardinality of a generating set of~$G$ (as a group). 

Bogopolski and Bux proved recently an analogue of Takahasi's Theorem
for fundamental groups of closed compact surfaces 
\cite[Proposition~2.2]{BB}. 
In~\cite{ASS}, Ara\'ujo, Silva and Sykiotis introduced the concept of
{\em Takahasi group}. A group $G$ is a Takahasi group if every 
ascending chain of subgroups of $G$ of bounded group rank is
stationary. In \cite{ASS}, among other results it is proved that:
\bt
\label{FinExtTakGp}
{\rm \cite[Theorem~4.1]{ASS}}
Every finite extension of a Takahasi group is a Takahasi group.
\et
It is also proved in \cite[Section~4]{ASS} that
every virtually free group is a Takahasi group, and 
the fundamental group of a finite graph of groups with
virtually polycyclic vertex groups and finite edge groups is a
Takahasi group.
On the other hand, in~\cite[Example~4.3]{ASS} it is shown that a free group 
of finite rank has arbitrarily long strict chains of subgroups 
of equal rank.

Takahasi's Theorem can be applied to prove that the subgroup 
$\per(\p)$ of periodic points of an automorphism $\p$ 
of a free group of finite rank is
finitely generated (which implies that the size of the periodic orbits
is bounded for each automorphism). 
Using the aforementioned generalization of Takahasi's 
Theorem, one obtains:

\bt
\label{tass}
{\rm \cite[Theorem~5.1]{ASS}}
Let $G$ be the fundamental group of a finite graph of groups with
finitely generated virtually nilpotent vertex groups and finite edge
groups. Then there exists a constant $M > 0$ such that
$$\rk_{\G}({\rm Per}(\p)) \leq M$$
for every $\p \in {\rm End}(G)$.
\et

As a consequence, a bound for the periods of a given endomorphism 
of~$G$ was also obtained in~\cite{ASS}.

In this paper, we consider the condition of Takahasi's Theorem in the
context of several varieties of semigroups, 
and apply results obtained to investigate the subsemigroup of fixed points  
and the subsemigroup of periodic points of an endomorphism of a semigroup 
of various kinds. 

The reader is assumed to have basic knowledge of semigroup theory,
universal algebra and automata theory.

We consider the following varieties throughout this paper:
\bi
\item
$\M$ -- the variety of all monoids (type (2,0)); 
\item
$\S$ -- the variety of all {\em semigroups} (type (2)); 
\item
$\I$ -- the variety of all {\em inverse semigroups} (type (2,1)); 
\item
$\G$ -- the variety of all {\em groups} (type (2,1)); 
\item
$\CR$ -- the variety of all {\em completely regular semigroups} (type
(2,1));
\item
$\C$ -- the variety of all {\em Clifford semigroups} (type
(2,1));
\item
$\CS$ -- the variety of all {\em completely simple semigroups} (type
(2,1)).
\ei 
The unary operation is $a \mapsto a\inv$ in the case of inverse
semigroups and groups, where $a\inv$ denotes the inverse of $a$, and
$a \mapsto \oo{a}$ in the case of completely regular semigroups, where
$\oo{a}$ is the unique inverse of $a$ commuting with it.   
Recall that $\CR$ contains both $\C$ and $\CS$. Also 
$$\G = \I \cap \CS, \quad \C = \I \cap \CR.$$

It is particularly important for us to remark which type of
subalgebras we have for each of these varieties: 
submonoids for~$\M$, subsemigroups for~$\S$, 
inverse subsemigroups for~$\I$, subgroups for~$\G$ and
completely regular subsemigroups for~$\CR$.
If~$\V$ is any of the varieties of type $\tau$ defined 
above and $S \in \V$, we write $T \leq_{\V} S$, 
and say that $T$ is a {\em $\V$-subalgebra\/} of~$S$, meaning  
that~$T$ is a $\tau$-subalgebra of $S$. 

\medskip


The paper is structured as follows.
In Section~\ref{Sec:TkSmgs}, we generalize the concept of Takahasi group to further
varieties of algebras. We show that finite $\J$-above semigroups are
Takahasi, and provide a full description of Takahasi completely simple
semigroups and of Takahasi Clifford semigroups. 
We also prove corollaries involving appropriate notions of finite index, 
as well as some negative results.

In Section~\ref{Sec:PeriodicPoints}, 
we introduce classes of semigroups $\UA$ (respectively $\UE$) where the rank of
fixed point subsemigroups is uniformly bounded for arbitrary
automorphisms (respectively endomorphisms). Using the results of Section~\ref{Sec:TkSmgs}, we
prove that the subsemigroups of periodic points are finitely generated
and periodic orbits are bounded for arbitrary endomorphisms
of finitely generated completely simple or Clifford semigroups whose
$\H$-classes are Takahasi groups and belong to~$\UE$. Similar results
are proved for balanced one-relator presentations of length 2.

\section{Takahasi semigroups} \label{Sec:TkSmgs}

Let $S$ be a semigroup and let $A \subseteq S$ be nonempty. 
We denote the subsemigroup of~$S$ generated by~$A$ by~$A^+$, and 
one has  
\[ 
A^+ = \bigcup_{n \geq 1} A^n\ .
\] 
If~$S$ is a monoid, the submonoid of~$S$ generated by~$A$ 
will be denoted by~$A^*$. Clearly, $A^* = A^+ \cup \{ 1 \}$.  
For any semigroup~$S$, the {\em rank\/} of $S$, 
denoted by~$\rk(S)$, is defined as 
$$\rk(S) = \min \bigl\{ |A| : \emptyset \neq A \subseteq
S,\; A^+ = S \bigr\},$$ 
if~$S$ is finitely generated, and as 
$\rk(S) = \infty$ otherwise.

Assume now that $\V$ is an arbitrary variety. Let
$S \in \V$ and $A \subseteq S$ be nonempty. We denote by 
$\langle A\rangle$ the $\V$-subalgebra of $S$ generated by $A$. 
Given $S \in \V$, we also define the {\em $\V$-rank\/} of $S$, 
denoted by~$\rk_\V (S)$, as 
$$\rk_{\V}(S) = \min \bigl\{ |A| : \emptyset \neq A \subseteq
S,\; \langle A\rangle = S \bigr\},$$ 
if~$S$ is finitely generated, and as 
$\rk_{\V}(S) = \infty$ otherwise.
Note that, 
for $\V \in \{\M,\I,\G,\CR,\C,\CS\}$, 
the inequalities
\beq
\label{inrank}
\rk_{\V}(S) \leq \rk(S) \leq 2\rk_{\V}(S)
\eeq
hold for every $S \in \V$. 
Thus Takahasi's Theorem could be
stated using the semigroup rank instead of the group rank. 

We generalize now the concept of Takahasi group for a variety
$\V$. Given $S \in \V$, we write $S \in \tak(\V)$ if every 
ascending chain of $\V$-subalgebras of $S$ of bounded $\V$-rank is
stationary. More precisely, whenever
$$T_{1}\leq T_{2}\leq \cdots \leq S$$
and $\rk_{\V}(T_n) \leq M$ for every $n \geq 1$, there exists some $p \geq
1$ such that $T_p = T_{p+1} = T_{p+2} = \cdots$

It follows easily from the definitions and (\ref{inrank}) that
$$\tak(\S) \cap \G \subseteq \tak(\G).$$
The following result shows that the opposite inclusion is far from
true. Given groups $G$ and $H$, we denote by $G \ast_{\G} H$ and $G
\ast_{\S} H$ the free product of $G$ and $H$ in the varieties $\G$ and
$\S$, respectively.

\bp
\label{notts}
\bi
\item[(i)] $\ZZ \times \ZZ \notin {\rm Tak}(\S)$;
\item[(ii)] if $H$ is a nontrivial group, then $\ZZ \ast_{\G} H \notin
  {\rm Tak}(\S)$; 
\item[(iii)] if $H$ is a nontrivial group, then $\ZZ \ast_{\S} H \notin
  {\rm Tak}(\S)$.
\ei
\ep

\bproof
(i) The group $\ZZ \times \ZZ$ is generated by $a = (1,0)$ and $b =
(0,1)$. For every $n \geq 1$,  
let $S_n$ be the subsemigroup of $\ZZ \times \ZZ$ generated by $a^{-2}$ and
$a^{2n-1}b$. Since $a^{2n-1}b = a^{-2}(a^{2n+1}b)$ for every
$n \geq 1$, then
\beq
\label{notts1}
S_1 \subseteq S_2 \subseteq S_3 \subseteq \dotsm
\eeq

Suppose that $a^{2n+1}b \in S_n$. Then the generator $a^{2n-1}b$ must
be used exactly once to get $a^{2n+1}b$, which is clearly impossible.
Hence $a^{2n+1}b
\notin S_n$ and so all the inclusions in (\ref{notts1}) are
strict. Therefore $\ZZ \times \ZZ
\notin {\rm Tak}(\S)$.

(ii) We use the same construction taking $a$ to be a generator of $\ZZ$
and $b \in H \setminus \{ 1 \}$. Once again, we have
(\ref{notts1}). 

Suppose that $a^{2n+1}b \in S_n$. Then
$$a^{2n+1}b = a^{-2k_0}a^{2n-1}ba^{-2k_1}a^{2n-1}b \dotsm
a^{-2k_{m-1}}a^{2n-1}ba^{-2k_m}$$
for some $m \geq 1$ and $k_0,\dotsc,k_m \geq 0$. Since there exists
always an odd number of $a$'s between any two consecutive $b$'s in the
right hand side, it follows easily from the normal form for 
the elements of the free product $\ZZ \ast_{\G} H$
that we must have $m = 1$, $k_1 = 0$ and $k_0 = -1$, a
contradiction. Hence $a^{2n+1}b
\notin S_n$ and so all the inclusions in (\ref{notts1}) are
strict. Therefore $\ZZ \ast_{\G} H
\notin {\rm Tak}(\S)$.

(iii) The proof of (ii) holds for $\ZZ \ast_{\S} H$ as well.
\eproof

Since $\ZZ \times \ZZ$ and the groups of the form $\ZZ \ast_{\G} H$ with
$H$ finite are Takahasi groups by \cite[Corollary~4.4]{ASS}, it
follows that $\tak(\G) \not\subseteq \tak(\S)$. By Proposition
\ref{notts}(ii), no free group of rank $> 1$ belongs to $\tak(\S)$. 

Another consequence of this last proposition 
is the bad behaviour of $\tak(\S)$ with respect to
direct products and free products. But first we discuss the case of
infinite cyclic groups.

\bp \label{Z-Takahasi}
The additive semigroup $\ZZ$ belongs to ${\rm Tak}(\S)$, 
and so does $(\N,+)$. 
\ep

For the proof we need some classical tools. 
Given a subsemigroup $S$ of the additive semigroup~$\N$ 
of natural numbers and $d \in \N$, 
we say that $S$ is {\em ultimately a $d$-segment\/} if 
there is $p \in \N$ such that for all $n \in \N$ such that $n \geq p$, 
we have $n \in S$ if and only if $d$ divides $n$. 
It is clear that $S$ cannot be ultimately a $d_1$-segment and 
ultimately a $d_2$-segment for distinct natural numbers $d_1$ and $d_2$.
Let 
\[
d_S = \min \bigl\{ \gcd\{x,y\} \colon  x, y \in S \bigr\}. 
\]
The next result can be found in~\cite[Sec.~II.4]{Gril} 
(see also \cite{JM, SS}).

\bt \label{SubsmgN}
If $S$ is a subsemigroup of $\N$, then 
$S \subseteq \N d_S$, $S$ is ultimately a $d_S$-segment, 
and $S$ is finitely generated.
\et

For a subsemigroup $S$ of $\N$, define 
\[
p_S = \min\bigl\{ p \in \N \colon \,
    \forall n \geq p \: ( d_S | n \Rightarrow n \in S)\bigr\} , 
\]
which is a natural number by Theorem~\ref{SubsmgN}. 

From Theorem~\ref{SubsmgN} and its dual for $\ZZ^-$ we 
can easily deduce the following corollary, which 
can also be found in~\cite[Sec.~II.4]{Gril}. 

\bc \label{subsmgsZ}
A subsemigroup of $\ZZ$ either contains only 
non-negative integers, or contains only 
non-positive integers, or is of the form $\ZZ d$ for some $d \in \N$, 
hence a subgroup of $\ZZ$. 
\ec

Now we are able to make the proof of Proposition~\ref{Z-Takahasi}. 

\begin{proof}[{\em\bf Proof of Proposition~\ref{Z-Takahasi}}]
First we prove that any infinite ascending chain of 
subsemigroups of~$\N$  
\[
S_1 \leq S_2 \leq \dotsm
\]
is stationary (this implies that $\N \in \tak(\S)$). 
Let us take such a chain. 
Then 
\[
d_{S_1} \geq d_{S_2} \geq \dotsm
\]
whence there is $k \in \N$ such that 
$d_{S_k} = d_{S_i}$ for all $i \geq k$. 
Given $i \geq k$, the fact that 
$S_i \leq S_{i+1}$ and $d_{S_i} = d_{S_{i+1}}$ 
implies that 
$p_{S_{i+1}} \leq p_{S_i}$. 
Thus 
\[
p_{S_k} \geq p_{S_{k+1}} \geq \dotsm , 
\]
and then there is 
$\ell \geq k$ such that $p_{S_\ell} = p_{S_i}$ 
for every $i \geq \ell$. 
Then in the chain 
\[
S_\ell \leq S_{\ell +1} \leq \dotsm
\]
any two of these semigroups only may differ 
in natural numbers less than $p_\ell$, and, hence, 
this chain is stationary.

Dually $\ZZ^{-}$ satisfies the same condition. 
Now let 
\beq \label{AscChainZ}
S_1 \leq S_2 \leq \dotsm
\eeq
be an infinite ascending chain of nontrivial subsemigroups of $\ZZ$. 
By Corollary~\ref{subsmgsZ}, either all these 
subsemigroups are contained in $\ZZ_0^-$, or all these 
subsemigroups are contained in $\N$, or there is 
$k \in \N$ such that $S_i$ is a subgroup of $\ZZ$ 
for every $i \geq k$.  
In the first two situations, the chain is stationary 
as we proved above.
In the third situation, the claim follows immediately from $\ZZ$ being
a noetherian ring.
Therefore $\ZZ \in \tak(\S)$. 
\end{proof}

As opposed to Proposition~\ref{Z-Takahasi}, we have the 
following. 

\bp \label{NotTak-Q}
The additive group $\Q$ of rational numbers is not 
a Takahasi group. 
\ep

\bproof
It suffices to observe that, defining,
for each positive integer $k$, the subgroup $H_k$ as 
being the cyclic subgroup of $\Q$ generated by  
$\frac{1}{2^k}$, we have the infinite ascending chain 
\[
H_1 < H_2 < H_3 < \dotsm 
\]
of subgroups of $\Q$ of rank~1. 
\eproof

A celebrated result of Group Theory, attributed to 
Higman, Neumann and Neumann,  and, independently 
to Freudenthall, states that every countable 
group is embeddable in a 2-generator group 
\cite{Gal,HNN,NN}. 
Thus, by Proposition~\ref{NotTak-Q}, there are 
finitely generated groups that are not Takahasi groups. 

Now we can prove:

\bp
\label{nop}
${\rm Tak}(\S)$ is not closed under:
\bi
\item[(i)] direct product;
\item[(ii)] free product.
\ei
\ep

\bproof
(i) This follows from Proposition \ref{notts}(i) and Proposition
\ref{Z-Takahasi}. 

(ii) Trivially, all finite semigroups belong to $\tak(\S)$. 
Now, the result follows from Proposition~\ref{notts}(iii) and 
Proposition~\ref{Z-Takahasi}.
\eproof

On the positive side, the following result provides a wide class of
examples of semigroups in $\tak(\S)$.

The quasi-order $\leq_{\J}$ on $S$ is defined
by
$$a \leq_{\J} b \hspace{.5cm} \mbox{if } a \in S^1bS^1.$$
A semigroup $S$ is {\em finite $\J$-above} if
$$\{ x \in S : x \geq_{\J} a \}$$
is finite for every $a \in S$.

\bt
\label{fjat}
Let $S$ be a finite $\J$-above semigroup. Then $S \in {\rm Tak}(\S)$.
\et

\bproof
Let $S$ be a finite $\J$-above semigroup 
and suppose that 
\beq 
\label{fjat1}
T_1 < T_2 < T_3 < \dotsm
\eeq
is an infinite ascending chain of finitely generated subsemigroups of $S$. 
It suffices to show that $\rk(T_n)$ is unbounded.

For each $n \geq 1$, we fix a generating set $A_n$ of $T_n$ of minimum size.
Let $a \in A_1$. Since $a \in T_n = A_n^+$ for every $n \geq 1$, we
have $a \in a_nT_n^1$ for some $a_n \in A_n$. Hence $a_n \geq_{\J} a$
in $S$ for every $n \geq 1$. Since $S$ is finite $\J$-above, it
follows that there exists a refinement 
$$T_{i_1} < T_{i_2} < T_{i_3} < \dotsm$$
of (\ref{fjat1}) and some 
$$b_1 \in A_{i_1} \cap A_{i_2} \cap A_{i_3} \cap \dotsm$$

Proceeding inductively, we assume now that there exists a refinement
\beq
\label{fjat2}
T_{j_1} < T_{j_2} < T_{j_3} < \dotsm
\eeq
of (\ref{fjat1}) and some distinct
$$b_1, \dotsc,b_n \in A_{j_1} \cap A_{j_2} \cap A_{j_3} \cap \dotsm$$
Since $T_{j_1} \subset T_{j_2}$, there exists some $c \in A_{j_2} \setminus
T_{j_1}$. 
Since $c \in T_{j_n} \setminus
T_{j_1} = A_{j_n}^+ \setminus
A_{j_1}^+$ for every $n \geq 2$, we
have $c \in T_{j_n}^1c_{n}T_{j_n}^1$ for some $c_n \in A_{j_n} \setminus
A_{j_1}$. Hence
$c_n \geq_{\J} c$ 
in $S$ for every $n \geq 2$. Since $S$ is finite $\J$-above, it
follows that there exists a refinement 
$$T_{k_1} < T_{k_2} < T_{k_3} < \dotsm$$
of (\ref{fjat2}) (and therefore of (\ref{fjat1})) and some 
$$b_{n+1} \in (A_{k_1} \cap A_{k_2} \cap A_{k_3} \cap \dotsm)\setminus
A_{j_1}.$$
Since $b_1, \dotsc,b_n \in A_{j_1}$, it follows that $b_1,
\dotsc,b_{n+1}$ are $n+1$ distinct elements of $A_{k_1} \cap A_{k_2}
\cap A_{k_3} \cap \dotsm\,$. By induction, such a property holds for
arbitrary $n$. In particular, for every $n \geq 1$, there exists some
$m \geq 1$ such that $|A_m| \geq n$. Thus $\rk(T_n) = |A_n|$ is
unbounded and so $S \in \tak(\S)$.
\eproof

A semigroup (monoid) presentation of the form $\langle A \mid u_i =
v_i\; (i \in I)\rangle$ is said to be {\em balanced} if
$|u_i| = |v_i|$ for every $i \in I$.

Since the semigroups in the statement of the next corollary are
clearly finite $\J$-above, we immediately get:

\bc
\label{ext}
The following semigroups belong to $\tak(\S)$:
\bi
\item[(i)] finite semigroups;
\item[(ii)] free semigroups and free monoids;
\item[(iii)] trace monoids;
\item[(iv)] semigroups or monoids defined by balanced presentations;
\item[(v)] free inverse semigroups and free inverse monoids.
\ei
\ec

Since $\tak(\S) \cap \I \subseteq \tak(\I)$, we get also:

\bc
\label{exti}
Free inverse semigroups and free inverse monoids belong to $\tak(\I)$.
\ec

We consider next $\CS$ and Rees matrix semigroups. But first
we need a lemma on ranks of groups defined by automata. 
Let $A$ be an alphabet. 
We denote by $A\inv$ a set of formal inverses of $A$. If
$M$ is a monoid of type $(2,1)$, and $x \mapsto x\inv$ is the unary
operation, we say that a monoid homomorphism $\p:(A \cup A\inv)^* \to
M$ is {\em matched} if $a\inv\p = (a\p)\inv$ holds for every $a \in
A$.

We say that $\A = (Q,q_0,T,E)$ is a {\em finite} $A$-{\em
automaton} if $Q$ is a finite set, $q_0 \in Q$, $T \subseteq Q$ and 
$E \subseteq Q \times A \times Q$, and  refer to the elements of $Q$ and
$E$ as {\em vertices} and {\em edges}, respectively.

We say that an $(A \cup A\inv)$-automaton $\A = (Q,q_0,T,E)$ is:
\bi
\item
{\em dual} if 
$$(p,a,q) \in E \iff (q,a\inv,p) \in E$$
holds for all $p,q \in Q$ and $a \in A$;
\item
{\em inverse} if it is dual, trim and deterministic;
\item
{\em Stallings} if it is inverse, $T = \{ q_0\}$ and the unique vertex which
may have outdegree 1 is $q_0$ \cite[Section~2]{BS1}.
\ei
Recall that the language recognized by $\A$ is the set 
$L(\A)$ of words $w \in (A \cup A^{-1})^*$ such that $w$ is the label of 
a path from $q_0$ to some $t \in T$. 
\bl
\label{ragr}
Let $A$ be a finite alphabet and let $\p:(A \cup A\inv)^* \to G$ be a
matched homomorphism onto a group. Let $\A = (Q,q_0,T,E)$ be a finite $(A \cup
A\inv)$-automaton such that $(L(\A))\p = G$. Then:
\bi
\item[(i)] ${\rm rk}_{\G}(G) \leq |E|$;
\item[(ii)] ${\rm rk}_{\G}(G) \leq |E| - |Q| + |\{ q_0\} \cup T|$ if
  $\A$ is trim. 
\ei
\el

\bproof
Since the trim part of $\A$ (i.e. the subautomaton induced by
all vertices lying in some successful path) has at most $|E|$ edges,
it is enough to consider the case when $\A$ is trim.

Let $\A_1 = (Q_1,q_0,q_0,E_1)$ be the automaton obtained by
identifying all the vertices of $T$ with~$q_0$. 
Clearly, $\A_1$ is trim.
Then let 
$\A_2 = (Q_1,q_0,q_0,E_2)$ 
be the automaton obtained from $\A_1$ by adding
edges of the form $p \mapright{a\inv} q$ whenever 
$E_1$ contains an edge $q \mapright{a} p$ $(a \in A\cup A\inv)$ but no
edge $p \mapright{a\inv} q$. Note that $\A_2$ is a trim dual automaton.

Next let $\A_3 = (Q_3,q_0,q_0,E_3)$ be the inverse automaton obtained
by successively identifying all pairs of edges of the form
$$q \mapleft{a} p \mapright{a} r$$
with $a \in A \cup A\inv$. This is the procedure known as {\em Stallings
foldings}. It is well known that the final result is independent of the
order in which foldings are executed (see \cite[Section~2]{BS1}).

Finally, let $\A_4 = (Q_4,q_0,q_0,E_4)$ be the Stallings automaton
obtained by successively removing from $\A_3$ all vertices 
of outdegree~1 that are distinct from $q_0$.

We prove that
\beq
\label{ragr1}
(L(\A_4))\p = (L(\A_3))\p = (L(\A_2))\p = (L(\A_1))\p = (L(\A))\p = G.
\eeq

We start with the equality $(L(\A_1))\p = (L(\A_0))\p$. 
Clearly, $L(\A) \subseteq L(\A_1)$, 
and so $(L(\A))\p \subseteq (L(\A_1))\p$.
To prove the opposite inclusion, it
suffices to assume that we are identifying $q_0$ with a single element
$t \in T$. 

We claim that
\beq
\label{ragr2}
\mbox{if $p \mapright{v} q$ is a path in $\A$ and $p,q \in \{
  q_0,t\}$, then } v\p \in G.
\eeq
Since $\A$ is trim, there exists some path $q_0 \mapright{u} t$ in
$\A$. If $p = t$ and $q = q_0$, then $uvu \in L(\A)$ and so
$v\p = (u\p)\inv(uvu)\p(u\p)\inv \in G$. The other cases are
straightforward variations of this one and can be omitted. Thus
(\ref{ragr2}) holds.

Now let $w \in L(\A_1)$. Then we may factor $w =
w_0w_1\dotsm w_n$ so that
$$q_0 \mapright{w_0} q_0 \mapright{w_1} \dotsm \mapright{w_n} q_0$$
enhances all the occurrences of the vertex $q_0$ in a path of $\A_1$
labelled by $w$. It follows that there are paths $p_i \mapright{w_i}
r_i$ in $\A$ with $p_i,r_i \in \{ q_0,t\}$ for $i = 0,\dotsc,n$. By
(\ref{ragr2}), we get $w_i\p \in G$ for every $i$, hence $w \in G =
(L(\A))\p$ and so $(L(\A_1))\p = (L(\A))\p$.

Similarly to the preceding equality, to prove the nontrivial inclusion
$(L(\A_2))\p \subseteq (L(\A_1))\p$, 
we may assume that $\A_2$ was obtained from $\A_1$ by adding the
single edge
$p \mapright{a\inv} q$. Let $w \in L(\A_2) \setminus L(\A_1)$. Then we
may factor $w = w_0a\inv w_1 \dotsm a\inv w_n$ so that
\beq
\label{ragr3}
q_0 \mapright{w_0} p \mapright{a\inv} q \mapright{w_1} \dotsm
\mapright{a\inv} q \mapright{w_n} q_0
\eeq
enhances all the occurrences of the new edge in a path of $\A_2$
labelled by $w$. Since $\A_1$ is trim, there exist paths of the form
$$q_0 \mapright{u} q, \quad p \mapright{v} q_0$$
in $\A_1$. Moreover, all the paths labelled by some $w_i$ in
(\ref{ragr3}) can be seen as paths in $\A_1$, hence $uav, w_0v, uw_n,
uw_iv \in L(\A_1)$ for $i = 1,\dotsc,n-1$ and we get
$$w\p = (w_0a\inv w_1\dotsm a\inv w_n)\p = (w_0v)\p((uav)\p)\inv
(uw_1v)\p\dotsm ((uav)\p)\inv (uw_n)\p \in G.$$
Thus $(L(\A_2))\p = (L(\A_1))\p$.

The equalities $(L(\A_4))\p = (L(\A_3))\p = (L(\A_2))\p$ are
well-known facts from the theory of Stallings foldings 
(\cite[Section~2]{BS1}), therefore (\ref{ragr1}) holds.

Let $FG_A$ denote the free group on $A$ and let $H$ be the subgroup of
$FG_A$ having Stallings automaton $\A_4$ (i.e. $H$ is the canonical
pre-image of $G$ in $FG_A$ with respect to $\p$). The famous rank formula for
Stallings automata states that 
$\rk_{\G}(H) = \frac{|E_4|}{2} -|Q_4|+1$, 
see \cite[Proposition~2.6]{BS1}). 
Since $G$ is a homomorphic image of $H$, it follows that 
\beq
\label{ragr4}
\rk_{\G}(G) \leq \frac{|E_4|}{2} -|Q_4|+1.
\eeq
Now each time we delete a vertex on constructing $\A_4$ from $\A_3$,
we delete at least two edges, hence
$$\frac{|E_4|}{2} -|Q_4|+1 \leq \frac{|E_3|}{2} -|Q_3|+1.$$
Similarly, each time we identify two vertices on constructing $\A_3$
from $\A_2$, we identify at least two pairs of edges, hence
$$\frac{|E_3|}{2} -|Q_3|+1 \leq \frac{|E_2|}{2} -|Q_1|+1.$$
Since $|E_2| \leq 2|E_1|$, condition~(\ref{ragr4}) and the above inequalities
yield
\beq
\label{ragr5}
\rk_{\G}(G) \leq |E_1| -|Q_1|+1.
\eeq
Now it suffices to note that $|E_1| \leq |E|$ and $|Q_1| = |Q| - |\{ q_0\}
\cup T| +1$.
\eproof

Let $S =
M[G,I,\Lambda,P]$ be a Rees matrix
semigroup. Given $X \subseteq S$, $i \in I$ and $\lambda \in
\Lambda$, we write
$$X^{(i\lambda)} = X \cap (\{ i \} \times G \times \{ \lambda\}).$$
Given $X \subseteq S$, we write
$$I_X = \bigl\{ i \in I : X \cap (\{ i \}
\times G \times \Lambda) \neq
\emptyset \bigr\}, \quad 
\Lambda_X = \bigl\{ \lambda \in \Lambda : X \cap (I
\times G \times \{ \lambda \}) \neq
\emptyset \bigr\}.$$

\bl
\label{rees}
Let $S = M[G,I,\Lambda,P]$ be a Rees matrix
semigroup. Let $T \leq_{\CS} S$, $i \in I_T$ and $\lambda \in
\Lambda_T$. Then:
\bi
\item[(i)] $T^{(i\lambda)}$ is a subgroup of $T$ and 
 $T^{(i\lambda)} \cong G^{(i\lambda)}$ for some $G^{(i\lambda)} \leq_{\G} G$;
\item[(ii)] if $T$ is finitely generated, 
then 
\beq
\label{rees4}
{\rm rk}_{\G}(T^{(i\lambda)}) \leq ({\rm rk}_{\CS}(T))^2+1.
\eeq
\ei
\el

\bproof
(i) Since $T \leq_{\CR} S$ makes
$T$ a union of groups and $\{ i \} \times G \times \{ \lambda\}$,
being an $\H$-class of $S$, is a
group itself, then $T^{(i\lambda)}$ is a subgroup of $T$.

Write $P = (p_{\lambda i})$. We define a mapping
$$\begin{array}{rcl}
\p:T^{(i\lambda)}&\to&G\\
(i,g,\lambda)&\mapsto&gp_{\lambda i}
\end{array}$$
For all $g,h \in G$, we have
$$((i,g,\lambda)(i,h,\lambda))\p = (i,gp_{\lambda i}h,\lambda)\p =
gp_{\lambda i}hp_{\lambda i} = ((i,g,\lambda)\p((i,h,\lambda)\p),$$
hence $\p$ is a group homomorphism. 

Since $\p$ is clearly injective, we get $T^{(i\lambda)}
\cong T^{(i\lambda)}\p$, so we may take $G^{(i\lambda)} =
T^{(i\lambda)}\p$.

(ii) Let $A$ be a generating set for $T$ (as a completely simple
semigroup) of minimum size. Note that
\beq
\label{igs}
I_T = I_A, \quad \Lambda_T = \Lambda_A.
\eeq

We take two new elements $q_0,t \notin \Lambda_A$ and let
$Q = \{ q_0,t\} \cup \Lambda_A$. Let also 
\begin{align*}
E &= \bigl\{ (q_0,g,\lambda') : (i,g,\lambda') \in A \bigr\}  
     \cup \bigl\{ (\lambda',p_{\lambda'i'}g,\lambda'') : 
           \lambda' \in \Lambda_A, \:(i',g,\lambda'') \in A \bigr\}
     \cup \bigl\{ (\lambda,p_{\lambda i},t) \bigr\} \\
  &\subseteq Q \times G \times Q.
\end{align*}
Then  $\A = (Q,q_0,t,E)$ is a {\em finite} $G$-{\em
  automaton}. Notions such as (successful) path and language
generalize from classical automata theory to $G$-automata in the
obvious way.

We prove that 
\beq
\label{rees1}
L(\A) = G^{(i\lambda)}.
\eeq

Let $g \in L(\A)$. Then there exists a path
\beq
\label{rees2}
q_0 \mapright{g_0} \lambda_0 \vvlongmapright{p_{\lambda_0i_1}g_1}
\lambda_1 \vvlongmapright{p_{\lambda_1i_2}g_2} \dotsm
\vvvvlongmapright{p_{\lambda_{n-1}i_n}g_n} \lambda_n = \lambda
\longmapright{p_{\lambda i}} t
\eeq
in $\A$, where $(i_j,g_j,\lambda_j) \in A$ for $j = 0,\dotsc,n$,
  $i_0 = i$ and
\beq
\label{rees3}
g = g_0p_{\lambda_0i_1}g_1p_{\lambda_1i_2}g_2\dotsc
p_{\lambda_{n-1}i_n}g_np_{\lambda i}.
\eeq
Hence
$$(i,gp_{\lambda i}\inv,\lambda) = (i,g_0,\lambda_0) (i_1,g_1,\lambda_1) \dotsm
(i_n,g_n,\lambda_n) \in A^+ = T$$
and so $(i,gp_{\lambda i}\inv,\lambda) \in T^{(i\lambda)}$. Thus $g =
(i,gp_{\lambda i}\inv,\lambda)\p \in T^{(i\lambda)}\p =
G^{(i\lambda)}$ and so $L(\A) \subseteq G^{(i\lambda)}$.

Conversely, let $g \in G^{(i\lambda)}$. 
Then $(i,gp_{\lambda
  i}\inv,\lambda) \in T^{(i\lambda)}$, therefore
there exist $(i_0,g_0,\lambda_0), (i_1,g_1,\lambda_1), \dotsc,$
$(i_n,g_n,\lambda_n) 
\in A$ such that
$$(i_0,g_0,\lambda_0) (i_1,g_1,\lambda_1) \dotsm
(i_n,g_n,\lambda_n) = (i,gp_{\lambda i}\inv,\lambda).$$
It is straightforward to check that this implies the existence of a
path of the form (\ref{rees2}) in $\A$. Since (\ref{rees3})
holds as well, it follows that $g \in L(\A)$ and so (\ref{rees1}) holds.

View $E$ as a finite alphabet and let $E\inv$ be a set of formal
inverses of $E$. Let $\theta:(E\cup E\inv)^* \to G$ be the matched
homomorphism which associates to each $e \in E$ its label. Replacing
$G$ by $G' = {\rm Im}(\theta)$, 
we may assume that $\theta$ is surjective.

Let
$\B$ be the finite $E$-automaton obtained by replacing the label of
each edge $e$ in $\A$ by $e$ itself. In view of (\ref{rees1}), we have
$G^{(i\lambda)} = L(\A) = (L(\B))\theta$. To prove that $\A$ is trim,
we take $\lambda' \in \Lambda_A$. Then $A$ contains elements of the
form $(i,g_1,\lambda_1)$, $(i_2,g_2,\lambda')$ and
$(i_3,g_3,\lambda)$. It is easy to check that there exists a path 
$$q_0 \mapright{g_1} \lambda_1 \vvlongmapright{p_{\lambda_1i_2}g_2}
\lambda' \vvlongmapright{p_{\lambda'i_3}g_3} \lambda
\longmapright{p_{\lambda i}} t$$
in $\A$, hence $\A$ is trim, and so is $\B$.

Thus we may apply Lemma \ref{ragr} to get
$${\rm rk}_{\G}(G^{i\lambda}) \leq |E| - |Q| + 2 \leq
|A|+|\Lambda_A|\cdot |A| +1 -|\Lambda_A| = (|\Lambda_A|+1)(|A|-1)+2.$$
Since $|\Lambda_A| \leq |A| = \rk_{\CS}(T)$, we get in view of (i)
$${\rm rk}_{\G}(T^{i\lambda}) \leq (\rk_{\CS}(T))^2+1$$ as required.
\eproof

\bt
\label{cst}
Let $S = M[G,I,\Lambda,P]$ be a Rees matrix
semigroup. 
Then, 
${S \in {\rm Tak}(\CS)}$ if and only if $G \in {\rm Tak}(\G)$.
\et

\bproof 
If $H$ is a subgroup of $G$, then $H \leq_{\CS} S$ up to isomorphism 
and $\rk_{\CS}(H) = \rk_{\G}(H)$. 
It follows that if
${S \in {\rm Tak}(\CS)}$, then $G \in {\rm Tak}(\G)$. 

Conversely, assume that $G$ is a Takahasi group. 
Let $N \geq 0$ and suppose 
that 
\beq
\label{cst1}
T_1 \leq T_2 \leq T_3 \leq \dotsm
\eeq
is an infinite chain of completely simple subsemigroups of $S$ with
$\rk_{\CS}(T_n) \leq N$ for every $n \geq 1$. If $\rk_{\CS}(T_n)$ is
realized by $A_n$, it follows from
(\ref{igs}) that
$$|I_{T_n}| = |I_{A_n}|\leq |A_n| = \rk_{\CS}(T_n) \leq
N.$$
Since 
$I_{T_1} \subseteq I_{T_2} \subseteq \dotsm$, then this chain must be
stationary. Similarly, the chain 
$\Lambda_{T_1} \subseteq \Lambda_{T_2} \subseteq \dotsm$ is
stationary. Removing finitely many terms of (\ref{cst1}) if needed, we
may assume that $I_{T_n} = I'$ and $\Lambda_{T_n} = \Lambda'$ for all
$n \geq 1$, for some $I'$ and $\Lambda'$ finite. 

In view of Lemma \ref{rees}(i), for all $i \in I'$ and $\lambda \in
\Lambda'$, we have chains of subgroups 
\beq
\label{cst2}
T_1^{(i\lambda)} \leq T_2^{(i\lambda)} \leq \dotsm
\eeq
By the proof of Lemma \ref{rees}(i), we get a chain
\beq
\label{cst3}
G_1^{(i\lambda)} \leq G_2^{(i\lambda)} \leq \dotsm
\eeq
of subgroups of $G$. 
Since Lemma \ref{rees}(ii) yields 
$$\rk_{\G}(G_n^{(i\lambda)}) =
\rk_{\G}(T_n^{(i\lambda)}) \leq N^2+1$$
and $G \in \tak(\G)$, each of the chains (\ref{cst3}) (and so each of
the chains (\ref{cst2})) must be stationary. Since $I'$ and $\Lambda'$
are both finite and 
$$T_n = \bigcup_{i \in I'} \bigcup_{\lambda \in \Lambda'}
T_n^{(i\lambda)}$$
for every $n \geq 1$, it follows that the chain (\ref{cst1}) is also
stationary. Therefore $S \in {\rm Tak}(\CS)$.
\eproof

\bc
Let $S$ be a completely simple semigroup 
and let $T \leq_{\CS} S$. Let $a \in T$ 
and let $G$ and $H$ be the $\H$-classes of $a$ 
in $S$ and in $T$, respectively. 
If $T \in {\rm Tak}(\CS)$ and 
$H$ is a subgroup of~$G$ of finite index, then 
$S \in {\rm Tak}(\CS)$. 
\ec

\bproof 
By Theorem~\ref{FinExtTakGp}, $H \in {\rm Tak}(\G)$ implies $G \in {\rm
  Tak}(\G)$. Now the claim follows from Theorem~\ref{cst}. 
\eproof

We consider next Clifford semigroups. This class of semigroups admits
various different characterizations. One of them states that a
semigroup $S$ is a Clifford semigroup if and only if $\H$ is a
semilattice congruence on $S$.

\bt
\label{fug}
Let $S$ be a Clifford semigroup. 
Then, $S \in {\rm Tak}(\C)$ if and only if 
every $\H$-class of~$S$ is a Takahasi group.
\et

\bproof
It is clear that  
if $S \in {\rm Tak}(\C)$, then
every $\H$-class of~$S$ is a Takahasi group.

Conversely, assume that every $\H$-class of $S$ 
is a Takahasi group. 
Let $H$ be an $\H$-class of $S$. 
First, we show that given $T \leq_{\C} S$ such that 
$T \cap H \neq \emptyset$, 
\beq
\label{fug1}
\rk_{\G}(T \cap H) \leq \rk_{\C}(T).
\eeq

Let $e$ denote the identity element of $H$. 
Since 
$T \cap H \neq \emptyset$ and $T$ is a Clifford subsemigroup of~$S$, 
then 
$T \cap H$ is a group of $H$.
Let $T' = \{ t \in T : te \in H \}$. 
We claim that
\beq
\label{fug3}
t \in T' \iff te \not<_{\J} e,
\eeq
where $<_{\J}$ denotes the $\J$-order in $T$. 

Indeed, we have always $te \leq_{\J} e$. Since $\J\, =\, \H$ in a
Clifford semigroup, we get
$$te \not<_{\J} e \iff te\,\J\, e \iff te\,\H\, e
\iff te \in H \iff t \in T'$$
and so (\ref{fug3}) holds.

If $t,u \in T'$, then
\beq
\label{fug2}
tue = t(eue) = (te)(ue) \in H,
\eeq
hence $T'$ is a subsemigroup of $T$. 

Assume first that $T' = T$. Then
$$\begin{array}{rcl}
\psi:T&\to&T\cap H\\
t&\mapsto&te
\end{array}$$
is a semigroup homomorphism by (\ref{fug2}). Given $t \in T$, and
since idempotents are central in a Clifford semigroup, we have
$$(t\psi)\inv = (te)\inv = e\inv t\inv = et\inv = t\inv e =
t\inv\psi,$$ hence $\psi$ is a homomorphism of Clifford
semigroups. Since $\psi$ fixes each
element of $T\cap H$, then it is surjective and so 
$$\rk_{\G}(T\cap H) = \rk_{\C}(T\cap H) \leq \rk_{\C}(T).$$

Thus we may assume that $T \setminus T' \neq \emptyset$. Let $(T\cap
H)^0$ be the Clifford semigroup obtained by adjoining a
zero element 0 to $T\cap H$. We define a mapping $\psi:T \to (T\cap
H)^0$ by 
$$t\psi = \left\{
\begin{array}{ll}
te&\mbox{ if }t \in T'\\
0&\mbox{ otherwise}
\end{array}
\right.$$
Let $t,u \in T$. If $t,u \in T'$, then 
$(tu)\psi = tue = (t\psi)(u\psi)$ by~(\ref{fug2}). 
Suppose next that $u \notin T'$. 
By~(\ref{fug3}), we have
$ue <_{\J} e$, hence $tue <_{\J} e$. 
It follows that $tu \not\in T'$, 
whence $(tu)\psi = 0 = (t\psi)(u\psi)$. 
Finally, assume that $t \notin T'$. 
Then $te <_{\J} e$ by (\ref{fug3}), hence $tue = teu <_{\J} e$ and so
$(tu)\psi = 0 = (t\psi)(u\psi)$. 
Thus $\psi$ is a semigroup homomorphism. 
Similarly to the case $T' = T$, we show that 
$(t\psi)\inv = t\inv\psi,$ hence $\psi$ is a homomorphism of Clifford
semigroups. Since $T' \neq T$ and $\psi$ fixes each
element of $T\cap H$, then it is surjective and so 
$$\rk_{\C}((T\cap H)^0) \leq \rk_{\C}(T).$$
If $A$ is a generating set of minimum size for $(T\cap H)^0$ in
$\C$, then $A \setminus \{ 0 \}$ generates $T\cap H$. 
Therefore~(\ref{fug1}) holds.

Now let $N \geq 0$ and suppose that
\beq
\label{fug4}
T_1 \leq T_2 \leq T_3 \leq \dotsm
\eeq
is an infinite chain of Clifford subsemigroups of $S$ with
$\rk_{\C}(T_n) \leq N$ for every $n \geq 1$. Consider the canonical homomorphism
$\p:S \to S/\H$ and write $Y_n = T_n\p$. 
Since the free semilattice on a set with $m$ elements has
$2^m-1$ elements, it follows from $\rk_{\C}(T_n) \leq N$ that 
$$|Y_n| \leq 2^{N}-1.$$
Therefore the chain
$Y' = Y_{1} \subseteq Y_{2} \subseteq \dotsm$ 
must be stationary. 
Removing finitely many terms of (\ref{fug4}) if needed, we
may assume that $Y_{n} = Y'$ for every $n \geq 1$. 
Thus $Y'$ consists of finitely many $\H$-classes 
$H_1,\dotsc,H_m$ of $S$ with $m \leq 2^N-1$.

For each $i = 1,\dotsc,m$, we get a chain of subgroups of $H_i$ of the form
\beq
\label{fug5}
T_1 \cap H_i \leq T_2 \cap H_i \leq \dotsm
\eeq
and from (\ref{fug1}) we have 
$\rk_{\G}(T_n \cap H_i) \leq \rk_{\C}(T_n) \leq N$ 
for every $n \geq 1$. 
Since $H_i \in \tak(\G)$, the chain (\ref{fug5}) must be
stationary for each $i$. 
As $Y'$ is finite and 
$$T_n = \bigcup_{i=1}^m (T_n \cap H_i)$$
for every $n \geq 1$, it follows that the chain (\ref{fug4}) is also
stationary. Therefore $S \in {\rm Tak}(\C)$.
\eproof

We have not succeeded so far on establishing whether a completely
regular semigroup where the $\H$-classes are Takahasi groups belongs
to $\tak(\CR)$. 
A first obstacle is that a finitely generated completely
regular semigroup may have infinitely many $\H$-classes. The first
such example was found by Clifford in \cite[Section~6]{Cli}: the free completely
regular semigroup on two generators.

\medskip

We introduce now a notion of index for Clifford semigroups. 
Let $S$ be a Clifford semigroup with $\H$-classes $\{ H_i : i \in I\}$
and let $T$ be a $(2,1)$-subalgebra of $S$. 
Then $T$ is also a Clifford semigroup and 
$T = \dotcup_{i \in I} (H_i \cap T)$. 
Thus, each $H_i \cap T$ is the empty set or a subgroup of $H_i$.
Define the {\em index of~$T$ in~$S$\/}, which we denote 
by $[S \colon\! T]$, by 
$$[S \colon\! T] = {\rm sup}\bigl\{ [H_i \colon\! H_i \cap T] \colon \,
i \in I\bigr\},$$  
with the convention that, for any group $G$, 
$[G \colon\! \emptyset]$ is the order of $G$. 
Clearly, this definition does not give rise to 
any contradiction if $S$ and $T$ are groups. 
Theorem~\ref{FinExtTakGp} can be generalized 
as follows. 

\bt \label{prop:Clif-SubsmgFiniteIndexTak}
If $S$ is a Clifford semigroup with a $(2,1)$-subalgebra $T$ 
of finite index such that $T \in {\rm Tak}(\C)$, then 
$S \in {\rm Tak}(\C)$.
\et

\bproof
Let $S$ and $T$ as in the statement. 
Write $S/\!\H\, = \{ H_i : i \in I\}$. 
Let $i \in I$. 
If $H_i \cap T = \emptyset$,
then $H_i$ is a finite group, by hypothesis, hence a 
Takahasi group. 
Besides, since $T \in \tak(\C)$, every 
nonempty group $H_i \cap T$ belongs to $\tak(\C)$, 
and therefore is a Takahasi group. 
Thus, by Theorem~\ref{FinExtTakGp}, 
$H_i$ is a Takahasi group.
Now Theorem~\ref{fug} gives the desired conclusion.
\eproof

Now we will compare this notion of index with a 
notion of index introduced by Gray and Ruskuc~\cite{GR08}.
Let $S$ be a semigroup and let $T$ be a subsemigroup of $S$. 
Define the binary relation $\L^T$ on $S$ by 
\[ 
a \mathrel{\L^T} b \,\iff\, T^1 a = T^1 b
\]
for all $a,b \in S$. 
Define $\R^T$ dually, and $\H^T = \L^T \cap R^T$. 
Each of these relations is an equivalence relation on $S$ 
and both $T$ and $S\setminus T$ are union of 
$\L^T$-classes (resp.\ $\R^T$-classes, $\H^T$-classes). 
In this context, those authors have defined the 
{\em Green index of $T$ in $S$\/}, which we denote by 
$[S \colon\! T]_\Gr$, as $n+1$, where $n$ is the cardinal of the set 
of $\H^T$-classes contained in $S\setminus T$.
This notion when restricted to groups $S$ and $T$ 
coincides with the usual notion of index of a subgroup in a group. 
Let us see how it relates with our notion of index 
in the case of Clifford semigroups. 

\bp \label{prop:Clif-Green-index}
Let $S$ be a Clifford semigroup such that 
$S/\!\H$ is finite and let $T$ be a 
$(2,1)$-subalgebra of $S$. 

If $[S \colon\! T] < \infty$, then 
$[S \colon\! T]_{\rm Gr} < \infty$.
\ep

\bproof
Write $S/\!\H\, = \{ H_i : i \in I\}$. Since $\H^T \subseteq \H$, each
$\H$-class $H_i$ is a union of 
$\H^T$-classes.

Suppose that $H_i \cap T \neq \emptyset$ and let $a, b \in H_i$. If 
$(H_i\cap T)a = (H_i \cap T)b$, 
then 
$a \in Tb$ and $b \in Ta$, 
whence 
$a \mathrel{\L^T} b$. 
Dually, if $a(H_i\cap T) = b(H_i \cap T)$,
then 
$a \mathrel{\R^T} b$.
Let 
\[
I_1 = \bigl\{ i \in I \colon\, H_i \cap T = \emptyset \bigr\} 
\qquad \text{ and } \qquad 
I_2 = I \setminus I_1.
\]
Let 
\[
\X = \bigcup_{i \in I_2} 
  \bigl\{ ((H_i\cap T)a, a(H_i \cap T)) \colon\, a \in H_i \bigr\}. 
\]
Then it is well defined the mapping 
$\psi \colon \X \to S/\!\H^T$ such that
$$((H_i \cap T)a, a(H_i\cap T)) \psi = H^T_a,$$ 
for $i \in I_2$ and $a \in H_i$, where $H^T_a$ denotes the $\H^T$-class of $a$. 

Assume that 
$[S \colon\! T] < \infty$. 
Then $H_i$ is finite for any $i\in I_1$. 
Moreover, $\X$ is finite, since $I$ is finite, and so is $[S \colon\! T]$, 
whence $\X\psi$ is finite too. 
Clearly  
\[ 
S/\!\H^T = 
   \Bigl( \,\bigcup_{i \in I_1} \bigl\{ H^T_a \colon\, 
   a \in H_i \bigr\} \Bigr) \cup  \X\psi , 
\] 
hence 
$[S \colon\! T]_\Gr \leq  \bigl\lvert S/\!\H^T \bigr\rvert < \infty$. 
\eproof

Next, we give an example that shows the analogue of 
Theorem~\ref{prop:Clif-SubsmgFiniteIndexTak} for 
Green index as well as  
the converse of Proposition~\ref{prop:Clif-Green-index} 
do not hold. 

\be
Let $G_0$ be a finitely generated group that is not 
a Takahasi group (we have observed that such a group exists). 
Let $A$ be a finite generating set of $G_0$ and let 
$G_1$ be the free group over~$A$. 
Then there is a surjective homomorphism $\phi \colon G_1 \to G_0$. 
Let $S = G_0 \,\dotcup\, G_1$ endowed with the product that 
extends the products in $G_0$ and in $G_1$ and such that, 
for $a \in G_0$ and $b \in G_1$, 
$a\cdot b = a (b\phi)$ and $b \cdot a = (b\phi) a$. 
Then $S$ is a Clifford semigroup (it is a strong semilattice of 
groups) such that 
$[S \colon\! G_1]_\Gr = 2$ and 
$[S \colon\! G_1] = |G_0| = \infty$. 
Moreover, by Theorem~\ref{taka}, $G_1$ is a Takahasi group. 
However, $S \not\in \tak(\C)$ since $G_0 \not\in \tak(\C)$.
\ee

\section{Periodic points} \label{Sec:PeriodicPoints}

In this section we apply the results of Section~\ref{Sec:TkSmgs} 
to the study of the subsemigroups of periodic points 
as well as of the periodic orbits of the endomorphisms of 
some classes of semigroups. 

For technical reasons, in this section we consider the empty set to be
a semigroup (of rank 0).

Let $\V$ be one of the varieties considered in Section~\ref{Sec:intro}. 
Given $S \in \V$ we denote by $\aut(S)$ 
(respectively $\endo(S)$) the automorphism group 
(respectively endomorphism monoid)
of~$S$. 
Note that, when dealing with homomorphisms, 
for the varieties of type~$(2,1)$ 
there is no need to specify the unary operation: 
any semigroup homomorphism between inverse (respectively
completely regular) semigroups preserves necessarily the respective
unary operation.

Given $\p \in \endo(S)$, its {\em fixed point} subsemigroup is
$$\fix(\p) = \{ a \in S : a\p = a\}$$
and its {\em periodic point} subsemigroup is
$$\per(\p) = \bigcup_{n\geq 1} \fix(\p^n).$$
Notice that $\fix(\p)$ and $\per(\p)$ 
are actually $\V$-subalgebras of~$S$.

Given $x \in \per(\p)$, the {\em period} of $x$ is the least $n \geq
1$ such that $x\p^n = x$.

Let $\UA(\V)$ denote the class of all $S \in \V$ such that
$$\exists N \in \N,\; \forall \p \in \aut(S),\hspace{.3cm} \rk(\fix(\p)) \leq N.$$
Similarly, we denote by $\UE(\V)$ the class of all $S \in \V$ such that
$$\exists N \in \N,\; \forall \p \in \endo(S),\hspace{.3cm}
\rk(\fix(\p)) \leq N.$$ 
Clearly, $\UE(\V) \subseteq \UA(\V)$. 
By considering the identity automorphism,
every $S \in \UA(\V)$ must be finitely generated. 
Note that, in view of~(\ref{inrank}), the definitions of 
$\UA(\V)$ and $\UE(\V)$ would not be affected if  
we had replaced $\rk$ by $\rk_{\V}$. 
In the case that~$\V$ is one of the varieties of type $(2,1)$, 
if $S \in \V$, then $\aut(S)$ (respectively $\endo(S)$) is 
formed by all semigroup automorphisms 
(respectively semigroup endomorphisms) of $S$, as we had observed. 
Thus, in this case we will refer to the classes $\UA(\V)$ and 
$\UE(\V)$ simply as $\UA$ and $\UE$, respectively. 
Observe, however, that semigroup homomorphisms between monoids do not necessarily respect 
the identity.

Let $FG_n$ denote the free group of rank $n \in \N$.
Using their sophisticated train track techniques, Bestvina and Handel
proved in \cite{BH} that, for every $\p \in \aut(FG_n)$, 
$$\rk(\fix(\p)) \leq n.$$
Latter, Imrich and Turner used this fact to prove in~\cite{IT} that 
the same relation holds for every $\p \in \endo(FG_n)$. 
Thus $FG_n \in \UE$ for every $n \in \N$. 

More generally, as stated in Theorem~\ref{tass}, fundamental groups of
finite graphs of groups with finitely generated virtually nilpotent
vertex groups and finite edge groups belong to $\UE$. 

For semigroups, we should mention that in the proof of 
\cite[Theorem~3.1]{RS} it was shown that, 
whenever $\p$ is an endomorphism of a
finitely generated {\em trace monoid} 
(i.e.\ partially commutative monoid) 
${\mathbb{M}}(A,I)$, 
we have 
$\rk_\M(\fix(\p)) \leq 2^{|A|}$. 
Therefore ${\mathbb{M}}(A,I) \in \UE(\M)$ if $A$ is finite.

Next, in a series of results, we provide some more instances of 
semigroups in~$\UE$.

\bl
\label{ue}
Let $S$ be a completely regular semigroup with finitely many
$\H$-classes. 
If all $\H$-classes of $S$ are in $\UE$, then $S \in \UE$.
\el

\bproof
Let $H_1, \dotsc, H_n$ be the $\H$-classes of $S$. For $i =
1,\dotsc,n$, assume that 
\beq
\label{ue3}
\rk(\fix(\psi)) \leq N_i \hspace{.5cm}\mbox{for every }\psi \in \endo(H_i).
\eeq
We show that
\beq
\label{ue1}
\rk(\fix(\p)) \leq \sum_{i = 1}^n N_i
\eeq
for every $\p \in \endo(S)$.

Fix $\p \in \endo(S)$ and let 
$$I = \bigl\{ i \in \{ 1,\dotsc,n \} : \fix(\p) \cap
H_i \neq \emptyset \bigr\}.$$ 
If $i \in I$, then $H_i\p \cap H_i \neq \emptyset$, and this 
yields $H_i\p \subseteq H_i$ since 
the $\H$-classes are the maximal subgroups of $S \in \CR$. 
For every $i \in I$, let $\p_{i}$ denote the
restriction of $\p$ to $H_i$, which is itself an endomorphism. It is
immediate that 
$$\fix(\p) = \bigcup_{i \in I} \fix(\p_{i}).$$
In view of (\ref{ue3}), we get
$$\rk(\fix(\p)) \leq \sum_{i \in I} \rk(\fix(\p_{i})) \leq \sum_{i \in
  I} N_i \leq \sum_{i = 1}^n N_i$$
and so (\ref{ue1}) holds. Therefore $S \in \UE$.
\eproof

\bp
\label{cue}
The following semigroups belong to $\UE$:
\bi
\item[(i)]
finitely generated completely simple semigroups with $\H$-classes in $\UE$;
\item[(ii)]
finitely generated Clifford semigroups with $\H$-classes in $\UE$.
\ei
\ep

\bproof
(i) Let $S = M[G,I,\Lambda,P]$ be finitely generated.  
Then both $I$ and $\Lambda$ must be finite. 
Thus $S$ has finitely many $\H$-classes and
the claim follows from Lemma~\ref{ue}.

(ii) Let $S \in \C$. 
Then the canonical mapping 
$S \to S/\!\H$  
is a surjective homomorphism, and $S/\!\H$ is a semilattice. 
Thus, if $S$ is finitely generated, then 
$S/\!\H$ is also finitely generated. 
Since finitely generated semilattices are well known to be finite, 
we may now apply Lemma \ref{ue}.
\eproof

It was noticed in~\cite{ASS} that there exist automorphisms 
$\p$ of the group $FG_2 \times \ZZ$ such that neither 
$\fix(\p)$ nor $\per(\p)$ is finitely generated as a group.  
Hence $FG_2 \times \ZZ \not\in \UA$. 
Now we give an example of a finite $\J$-above semigroup 
which satisfies the analogous property. 

\be
\label{exth}
Let $S$ be the semigroup defined by the presentation 
$\langle\, a,b,c \mid cac = cbc \, \rangle$. 
This semigroup is finite $\J$-above, since 
$|cac| = |cbc|$. 
Clearly, there exists an endomorphism $\p$ of $S$ 
satisfying $a\p = b$, $b\p = a$ and $c\p = c$. 
Since $\p^2 = \text{\rm id}_S$, 
this homomorphism is an automorphism of~$S$. 
The elements of $S$ 
\[ 
cac,\, (ca)^2 c,\, (ca)^3 c, \dotsc 
\]
are pairwise distinct and belong to ${\rm Fix}(\p)$. 
By definition of $S$, given $n \in \N$, the nontrivial 
factorizations of $(ca)^n c$ in~$S$ of length two are 
$(ca)^n c = uv$, where 
\[ 
(u,v) = \bigl( (ca)^k c, (ac)^{n-k} \bigr) 
\quad \text{or} \quad  
(u,v) = \bigl( (ca)^{k+1} , (ca)^{n-k-1}c \bigr), 
\]
with $k \in \{ 0, \dotsc, n-1\}$. 
However, in any of these situations, 
$\{ u,v \} \not\subseteq {\rm Fix}(\p)$. 
Then any generating set of ${\rm Fix}(\p)$ 
contains $cac,\, (ca)^2 c,\, (ca)^3 c, \dotsc$, 
and hence ${\rm Fix}(\p)$ is not finitely generated. 
Therefore $S \not\in \UA$. 
Notice that ${\rm Per}(\p) = {\rm Fix}(\p)$, since $\p^2 = \p$. 
\ee

The following result and its corollary show that the above counterexample is in
some sense minimal among the semigroups not in 
$\UA$ defined by one-relator balanced presentations.

\bt
\label{ltwo}
Let $M$ be the monoid defined by a 
finite presentation of the form 
\beq 
\label{ltwo2}
\langle A \mid a_1a_2 = a_3a_4\rangle,
\eeq
with $a_1, \dotsc,a_4 \in A$ not necessarily
distinct. 
Let $\p \in {\rm End}(M)$. Then
$${\rm rk}_\M({\rm Fix}(\p)) \leq |A|.$$
\et

\bproof
We use induction on $|A|$. 
The case $|A| = 1$ is trivial, since $\fix(\p) = \{ 1\}$ or 
$\fix(\p) = M$. 
Now assume
that $|A| > 1$ and the claim holds for smaller alphabets.

The possibility of induction is legitimate since any submonoid of $M$
generated by a proper subset of $A$ can still be defined by a
presentation of the form (\ref{ltwo2}) as we prove next. 
Let $A'$ be such a
subset, and let $M'$ be the submonoid of $M$ generated by $A'$.
If $a_1, \dotsc,a_4 \in A'$, then it is easy to see that $M'$ is
presented by $\langle A' \mid a_1a_2 = a_3a_4\rangle$, and $M$ is the
free product of $M'$ and the free monoid on $A \setminus A'$. Assume now that
$\{ a_1, \dotsc,a_4 \} \not\subseteq A'$. 

Suppose that $a_1 \neq a_3$ and $a_2 \neq a_4$. If we view (\ref{ltwo2}) as
a group presentation, it becomes a one-relator presentation with a
cyclically reduced relator $a_1a_2a_4\inv a_3\inv$. By Magnus' famous
Freiheitssatz (see \cite{MKS}), since the subgroup
generated by $A'$ misses one of the generators occurring in the
cyclically reduced relator, it is the free group on $A'$. Now it
follows easily that $M'$ is the free
monoid on $A'$, hence trivially definable by a presentation of the
form (\ref{ltwo2}). 

Finally, by left-right symmetry, we only need to consider the case of
presentations of the form $\langle A \mid ab = a^2\rangle$ or $\langle
A \mid ab = ac\rangle$. In the first case, we may still use the
Freiheitssatz since $ba\inv$ is a cyclically reduced relator where $a$
and $b$ both occur. The second case follows easily from the fact
that there are no nontrivial overlappings involving $ab$ and $ac$,
thus every application of the relation $ab = ac$ (involving a letter
which is not in $A'$) must be ``undone'' the exact same way. That is, 
$M'$ is the free monoid on $A'$.

Let $A_f$ be the set of letters of $A$ occurring in any word
representing any fixed point of $\p$. 
Let $M'$ denote the submonoid of $M$ generated by $A_f$. 
Then $u = u\p$ for all $u \in \fix(\p)$ yields
$A_f\p \subseteq M'$, and so the restriction $\p' = \p|_{M'}$ is an
endomorphism of $M'$. Moreover, $\fix(\p) = \fix(\p')$. 

If $A_f \subset A$, we get 
$$\rk_\M(\fix(\p)) = \rk_\M(\fix(\p')) \leq |A_f| < |A|$$
by the induction hypothesis. Now assume that $A_f = A$.

Before proceeding, let us notice that, 
since words representing the same element of $M$ must have the same
length, 
we have a natural concept of {\em length} for the elements of $M$.

Then, given $a \in A$, there exist $x, y \in M$ such that 
$xay \in \fix(\p)$.  
Thus, since $xay = (xay)\p^n$ for all $n \geq 1$, 
the element $a\p^n$ is a factor of $xay$, and 
$|a\p^n| \leq |xay|$ for all $n \geq 1$.
Therefore 
$\{ a,a\p, a\p^2, \dotsc \}$ must be finite.


It follows that, for
every $a \in A$, there exist $m_a \geq 0$ and $p_a \geq 1$ such that
$a\p^{m_a+p_a} = a\p^{m_a}$. 
For any integers $m$ and $p$ such that 
$m \geq m_a$ and $p$ is positive multiple of $p_a$ for all $a \in A$, 
we have 
$a\p^{m+p} = a\p^{m}$ for every $a \in A$. 
Thus, we may take such an $m$ and such a $p$ satisfying 
$m = p-1$, yielding 
$a\p^{2p-1} = a\p^{p-1}$ for all $a \in A$. 
Hence, for every $u \in M$,
\beq
\label{ltwo1}
u\p^{2p-1} = u\p^{p-1}
\eeq
and, in general, $u \p^{jp-1} = u \p^{p-1}$ for 
every $j \geq 1$. 
It follows that 
\beq
\label{ltwo4}
u\p^n = 1 \Rw u\p^{p-1} = 1
\eeq
for all $u \in M$ and $n \geq 1$.
We now prove that
\beq
\label{ltwo3}
\fix(\p) = (\fix(\p^{p-1}))\p^{p}.
\eeq

If $u \in \fix(\p)$, then $u = u\p^{p-1} = u\p^p$ yields $u \in
(\fix(\p^{p-1}))\p^{p}$. 
Conversely, let $v \in \fix(\p^{p-1})$ and $u = v\p^p$. 
Then, in view of~\eqref{ltwo1}, 
$u\p = v\p^{p+1} = v\p^{p-1}\p^{p+1} = v\p^{2p} = v\p^{p} = u$, 
and so (\ref{ltwo3}) holds.

Thus $\rk_\M(\fix(\p))\leq \rk_\M(\fix(\p^{p-1}))$. 
By replacing $\p$ by $\p^{p-1}$, in view of~\eqref{ltwo4} 
this allows us to assume that 
$$a\p^n = 1 \Rw a\p = 1$$
for all $a \in A$ and $n \geq 1$. 
Let
$$A_0 = A \cap 1\p\inv 
\qquad \text{ and } \qquad 
A_1 = A\setminus A_0\, .$$
Then 
\beq
\label{ltwo5}
A_1\p^n \subseteq M \setminus A_0^*
\eeq
for every $n \geq 1$. 

Now we split our discussion into two cases. We consider first the case
$A_0 \neq \emptyset$. 

Consider the homomorphism between free monoids 
$\pi:A^* \to A_1^*$ which erases the letters of~$A_0$.
Let $M_1$ be the monoid defined by the presentation
\beq
\label{simpl2}
\langle A_1 \mid (a_1a_2)\pi= (a_3a_4)\pi\rangle.
\eeq
We claim that (\ref{simpl2}) is equivalent to some presentation of the
form (\ref{ltwo2}). 
This clearly holds if 
$|(a_1a_2)\pi| = |(a_3a_4)\pi|$. 
On the other hand, the facts that   
$1\pi^{-1} = A_0^*$, $(a_1a_2)\p = (a_3a_4)\p$ and 
in~$M$ there is no invertible 
elements other than~$1$ imply
that $(a_1a_2)\pi = 1$ if and only if $(a_3a_4)\pi = 1$. 
Therefore we are left, in view of left-right symmetry, with the case 
$a_1,a_3,a_4 \in A_1$ and $a_2 \in A_0$.

Suppose that $a_1 \in \{ a_3,a_4\}$. Then $|(a_1a_2)\p| =
|(a_3a_4)\p|$ implies $a_i\p = 1$ for some $i \in \{ 3,4\}$,
contradicting $a_i \in A_1$. Hence $a_1 \notin \{ a_3,a_4\}$. 
But then $M_1$ is the free monoid on
$A_1 \setminus \{ a_1\}$. 

Let $\theta:A^* \to M$ and $\theta_1:A_1^* \to M_1$ be the canonical
homomorphisms. 
Since both homomorphisms 
$\pi{\theta_|}_{\scriptsize A_1^*} \p$ and $\theta\p$ 
coincide for letters of 
$A_0$ and $A_1$, we have
\beq
\label{simpl5}
\pi{\theta_|}_{\scriptsize A_1^*} \p = \theta\p.
\eeq
As $\ker(\theta)$ is the congruence generated by the
relation $a_1a_2 = a_3a_4$ and $(a_1a_2)\pi\theta_1 =
(a_3a_4)\pi\theta_1$, there exists a homomorphism $\pi':M \to M_1$
such that 
\beq
\label{simpl3}
\theta\pi' = \pi\theta_1.
\eeq
On the other hand, since $\ker(\theta_1)$ is the congruence generated by the
relation 
$(a_1a_2)\pi = (a_3a_4)\pi$ 
and 
$(a_1a_2)\pi{\theta_|}_{\scriptsize A_1^*} \p = 
 (a_3a_4)\pi{\theta_|}_{\scriptsize A_1^*} \p$ 
in view of (\ref{simpl5}), 
there exists a homomorphism 
$\psi:M_1 \to M$ 
such that 
\beq
\label{simpl4}
\theta_1\psi = {\theta_|}_{\scriptsize A_1^*} \p.
\eeq
We show that
\beq
\label{simpl6}
\pi'\psi = \p.
\eeq

Indeed, since $\theta$ is onto, (\ref{simpl6}) follows from
$$\theta\pi'\psi = \pi\theta_1\psi = \pi{\theta_|}_{\scriptsize A_1^*} \p = \theta\p,$$
where these equalities come from 
(\ref{simpl3}), (\ref{simpl4}) and (\ref{simpl5}), respectively. 

We show next that
\beq
\label{simpl7}
\fix(\p) = (\fix(\psi\pi'))\psi.
\eeq

Let $v \in A^*$ be such that $v\theta \in \fix(\p)$. By
(\ref{simpl5}) and (\ref{simpl4}), we have
$$v\theta = v\theta\p = v\pi\theta\p = v\pi\theta_1\psi.$$
Now (\ref{simpl4}), (\ref{simpl5}) and (\ref{simpl3}) yield
$$(v\pi\theta_1)\psi\pi' = v\pi\theta\p\pi' = v\theta\p\pi' =
v\theta\pi' = v\pi\theta_1,$$
hence 
$v\pi\theta_1 \in \fix(\psi\pi')$ 
and 
$\fix(\p) \subseteq (\fix(\psi\pi'))\psi.$ 

Conversely, let $v \in A_1^*$ be such that $v\theta_1 \in
\fix(\psi\pi')$. Then (\ref{simpl6}) yields
$$v\theta_1\psi\p = v\theta_1\psi\pi'\psi = v\theta_1\psi,$$
hence $(\fix(\psi\pi'))\psi \subseteq \fix(\p)$ and (\ref{simpl7})
holds. 

Now we may apply the induction hypothesis to 
the endomorphism $\psi\pi'$ of $M_1$ 
to get
$$\rk_\M(\fix(\psi\pi')) \leq |A_1|.$$
Therefore (\ref{simpl7}) yields
$$\rk_\M(\fix(\p)) \leq \rk_\M(\fix(\psi\pi')) \leq |A_1| < |A|$$
and the case $A_0 \neq \emptyset$ is settled.

We assume now that $A_0 = \emptyset$. 
By (\ref{ltwo5}), we have $|a\p| \geq 1$ for every $a \in A = A_1$. 
Recall that we are considering the case $A = A_f$. 
Thus, if there exists $a \in A$ such that $|a\p| > 1$, then 
$|u\p| > |u|$ for every $u \in M$ that has $a$ as a factor, 
contradicting the fact that $a \in A_f$. 
It follows that 
\beq
\label{simpl1}
|a\p| = 1 \mbox{ for every }a \in A.
\eeq 


Let $B = A \setminus \{a_1,\dotsc ,a_4\}$ 
and 
$C = \{a_1,\dotsc ,a_4\}$. 
Any word $u$ of $A^*$ can be factorized in a unique way 
in the form 
$u = w_0 u_1 w_1 \dotsm u_n w_n$, 
where 
$n \geq 0$, $w_0, w_n \in B^*$, $w_1, \dotsc, w_{n-1} \in B^+$ and  
$u_1, \dotsc, u_n \in C^*$. 
A word $v$ of $A^*$ represents the same element of~$M$ as such a word $u$ 
if and only if 
$v = w_0 v_1 w_1 \dotsm v_n w_n$, 
where 
$v_1, \dotsc, v_n \in C^*$ 
are such that $u_i$ and $v_i$ represent the same element of~$M$ 
for every $i = 1, \dotsc, n$.
Then, in view of~\eqref{simpl1} and that 
every letter of~$B$ occurs in a fixed point of~$\p$, 
we have 
\beq
\label{ltwo11}
a\p = a \mbox{ for every }a \in B, 
\eeq 
and $C\p \subseteq C$. 
Let $M_{C\p}$ be the submonoid of $M$ generated by $C\p$. 
Then 
\[ 
\fix(\p) = \bigl( B \cup \fix\bigl({\p_|}_{M_{C\p}}\bigr) \bigr)^* . 
\] 

If $C\p \neq C$, by the induction hypothesis it follows that 
\[
\rk_\M(\fix(\p)) = |B| + \rk_\M\bigl( \fix\bigl({\p_|}_{M_{C\p}}\bigr) \bigr)  
  \leq |B| +|C\p| < |B| + |C| = |A|.
\]

Suppose now that $C\p = C$. 
Then $\p$ induces a permutation on~$A$. 
If $\p$ is the identity of~$M$, then 
$\fix(\p) = M$, whence $\rk_\M(\fix(\p)) = |A|$. 
If $\p$ is not the identity of~$M$ and the words 
$a_1 a_2$ and $a_3 a_4$ are equal, then 
$M$ is the free monoid on~$A$, and 
$\fix(\p) = B^*$, whence 
$\rk_\M(\fix(\p)) = |B| < |A|$.
We proceed under the assumption that $\p$ is not 
the identity of~$M$ and 
that the words $a_1 a_2$ and $a_3 a_4$ are distinct. 
In view of (\ref{ltwo11}), and using left-right symmetry, we
may assume that $a_1\p \neq a_1$. 
Suppose that $a_1\p \neq a_3$. 
Then, from $(a_1a_2)\p = (a_3a_4)\p$ and the definition of~$M$, 
we have $a_1\p = a_3\p$ and $a_2\p = a_4\p$, whence 
$a_1 = a_3$ and $a_2 = a_4$, a contradiction.
Hence $a_1\p = a_3$. 
It follows that $a_1 \neq a_3$. From the fact that $\p$ is induces a permutation on~$A$ and 
the definition of~$M$, the homomorphism~$\p$ must permute 
$a_1$ with $a_3$ as well as $a_2$ with $a_4$. 
We split our discussion into two cases.

In the case that $a_2 = a_4$ it is easy to check that  
$$\fix(\p) = \bigl( (A\setminus \{ a_1,a_3\})\cup \{ a_1 a_2\} \bigr)^*.$$
Therefore $\rk_\M(\fix(\p)) < |A|$ in this case.

Assume now that $a_2 \neq a_4$. 
Then our defining relation of $M$ must be of the form 
$ab = cd$ or $ab = ba$ or $a^2 = b^2$, with $a,b,c,d$ distinct. 
If it is $ab = cd$ or $ab = ba$, it is easy to check that 
$$\fix(\p) = (B \cup \{ ab \})^*,$$ 
yielding $\rk_\M(\fix(\p)) < |A|$. 
Let us consider the case where the relation is $a^2 = b^2$. 
The rewriting system (see \cite{BO} for details)
$$\bigl\{ b^2 \longrightarrow a^2, \: ba^{2} \longrightarrow a^{2}b \bigr\} $$ 
is {\em noetherian} (there are no infinite chains of
reductions since the lexicographic order is a well-order)  
and {\em locally confluent} (since the unique overlappings between
relators are those of the form 
$bb^2 = b^2b$ and $b(ba^{2}) = b^2a^{2}$, 
and both lead to commutative diagrams such as
$$\xymatrix{ 
b^3 \ar[r] \ar[d] & ba^2 \ar[dl] & & b^2a^{2} \ar[d] \ar[r] & ba^{2}b
\ar[r] & a^{2}b^2 \ar[dll] \\
a^2b &&& a^{4} &&
}$$
The rewriting system is then {\em confluent} and the set of
irreducible words is a set of normal forms for~$M$. 
The irreducible words are 
those of the form $w_0 v_1 w_1 \dotsm v_n w_n$, 
where 
$n \geq 0$, $w_0, w_n \in B^*$, $w_1, \dotsc, w_{n-1} \in B^+$ and  
$v_1, \dotsc, v_n \in a^*(ba)^*\{ 1, b\}$.

Now, it is easy to check that 
$$\fix(\p) = \bigl( B \,\cup\, \bigr\{ a^2 \bigl\} \bigr)^*,$$ 
and $\rk_\M(\fix(\p)) = |B|+1 < |A|$ as required.
\eproof

Now, we may conclude the following.

\bc \label{Cor-ltwo}
Any monoid (respectively semigroup) defined by a 
finite presentation of the form 
\[
\langle A \mid a_1a_2 = a_3a_4\rangle,
\]
with $a_1, \dotsc,a_4 \in A$ not necessarily
distinct, is in $\UE(\M)$ (respectively $\UE(\M)$).  
\ec

\bproof
The statement for monoids follows directly from 
Theorem~\ref{ltwo}. 
Suppose that $S$ is a semigroup defined by  
such a (semigroup) presentation 
$\langle A \mid a_1a_2 = a_3a_4\rangle$. 
Since~$S$ does not have an identity, 
the monoid $S^1$ is also defined by the (monoid) 
presentation 
$\langle A \mid a_1a_2 = a_3a_4\rangle$. 
Then any (semigroup) endomorphism~$\p$ of~$S$ 
can be naturally extended to a (monoid) endomorphism~$\p_1$ of~$S^1$.  
For such endomorphisms, we have 
$\fix(\p_1) = \fix(\p) \cup \{ 1 \}$, 
whence
$\rk(\fix(\p)) = \rk_\M (\fix(\p_1))$. 
The desired conclusion now follows from Theorem~\ref{ltwo}. 
\eproof

Contrarily to what happens in free groups, 
Rodaro and the third author proved 
\cite[Theorem~3.10 and Corollary~3.11]{RS2}
that any nontrivial finitely generated free inverse monoid 
has automorphisms~$\p$ such that $\fix(\p)$ is not 
finitely generated. 
Hence nontrivial finitely generated free inverse monoids 
are not in $\UA(\M)$, and thus not in $\UE(\M)$ either. 

Next we will see some relationships between 
the classes $\tak(\V)$ and the 
classes $\UA$ and $\UE$.

A simple adaptation of an argument known for groups (see e.g.\ the
proof of \cite[Theorem~5.1]{ASS}) allows us to prove the following:

\bt
\label{fgp}
Let $\V$ be one of the varieties of type~$(2)$ or~$(2,1)$ 
considered in Section~\ref{Sec:intro} and let $S \in {\rm Tak}(\V)$.
\bi
\item[(i)] If $S \in \UA$, then ${\rm Per}(\p)$ is finitely generated
  for every $\p \in {\rm Aut}(S)$.
\item[(ii)] If $S \in \UE$, then ${\rm Per}(\p)$ is finitely generated
  for every $\p \in {\rm End}(S)$.
\ei
\et

\bproof
(i) If $S \in \UA$, then there exists some $N \in \N$ 
such that 
\beq
\label{fgp1}
\forall \psi \in \aut(S),\hspace{.3cm} \rk(\fix(\psi)) \leq N.
\eeq

Let $\p \in \aut(S)$. 
It is easy to see that
\beq 
\label{psv1} 
m|n
\hspace{.3cm} \Rw \hspace{.3cm} \fix(\p^{m}) \leq \fix(\p^{n}) 
\eeq
holds for all $m,n \geq 1$.
Hence we have an ascending chain of subsemigroups of $S$ of the form
$$\fix(\p) \leq \fix(\p^{2!}) \leq \fix(\p^{3!}) \leq \dotsm$$
By (\ref{fgp1}), we have $\rk(\fix(\p^{n!})) \leq N$ 
for every $n \geq 1$. 
Since $S \in \tak(\V)$,
there exists some $k \geq 1$ such that $\fix(\p^{n!}) =
\fix(\p^{k!})$ for every $n \geq k$. In view of (\ref{psv1}), we get
$$\per(\p) = \bigcup_{n \geq 1} \fix(\p^{n}) = \bigcup_{n \geq 1}
\fix(\p^{n!}) = \fix(\p^{k!}).$$
Therefore $\rk(\per(\p)) = \rk(\fix(\p^{k!})) \leq N$ by (\ref{fgp1})
and so $\per(\p)$ is finitely generated.

(ii) Similar.
\eproof

We remark that, even for $S \in \tak(\V)$, the conditions $S \in \UA$
or $S \in \UE$ are far from necessary to get finitely generated
periodic subalgebras. For instance, it follows from the results in
\cite[Section~3]{RS2} that 
$\fix(\psi)$ is not finitely generated when $M$ is the free monogenic
inverse monoid and $\psi \in \aut(M)$ sends the generator $a$ to its
inverse $a\inv$. However, by \cite[Theorem~3.8]{RS2}, $\per(\p)$ is finitely
generated for every endomorphism $\p$ of a free inverse monoid of
finite rank.

A straightforward adaptation of the proof of \cite[Corollary~5.2]{ASS}
yields the following corollary. We include the (short) proof for
completeness.  

\bc
\label{boup}
Let $\V$ be one of the varieties of type~$(2)$ or~$(2,1)$ 
considered in Section~\ref{Sec:intro} and let $S \in {\rm Tak}(\V)$.
\bi
\item[(i)] If $S \in \UA$ and $\p \in {\rm Aut}(S)$, then there
  exists a constant $R_{\p} > 0$ such that every $a \in {\rm Per}(\p)$
  has period less or equal to~$R_{\p}$. 
\item[(ii)] If $S \in \UE$ and $\p \in {\rm End}(S)$, then there
  exists a constant $R_{\p} > 0$ such that every $a \in {\rm Per}(\p)$
  has period less or equal to~$R_{\p}$. 
\ei
\ec

\bproof
(i) By Theorem \ref{fgp}, we may write 
$\per(\p) = \{ a_1, \dotsc, a_r\}^+$. 
Let $R_{\p}$ denote the least common multiple of the periods
of the elements $a_1, \dotsc, a_r$. 
Let $a \in \per(\p)$. 
Then there exist 
$i_1, \dotsc, i_n \in \{ 1, \dotsc, r\}$ 
such that 
$a = a_{i_1}\dotsm a_{i_n}$. 
It follows that
$$a\p^{R_{\p}} = (a_{i_1}\dotsm a_{i_n})\p^{R_{\p}} = (a_{i_1}\p^{R_{\p}})\dotsm
(a_{i_n}\p^{R_{\p}}) = a_{i_1}\dotsm a_{i_n} = a,$$
hence $a$ has period less or equal to~$R_{\p}$. 

(ii) Similar.
\eproof

We note that Theorem~\ref{fgp} and Corollary~\ref{boup} 
also hold for $\V = \M$ by replacing $\UA$ and $\UE$ by 
$\UA(\M)$ and $\UE(\M)$, respectively. 

Now we get the following result:

\bt
\label{perfh}
Let $S \in \CS \cup \C$ be finitely generated with all $\H$-classes in
${\rm Tak}(\G) \cap \UE$. Then:
\bi
\item[(i)]
${\rm Per}(\p)$ is finitely generated for every $\p \in {\rm End}(S)$;
\item[(ii)]
for every $\p \in {\rm End}(S)$, there
  exists a constant $R_{\p} > 0$ such that every $a \in {\rm Per}(\p)$
  has period less or equal to~$R_{\p}$. 
\ei
\et

\bproof
By Theorems~\ref{cst} or~\ref{fug}, we have accordingly $S \in
\tak(\CS)$ or $S \in \tak(\C)$. By Proposition~\ref{cue} and 
Theorem~\ref{fgp}(ii), 
${\rm Per}(\p)$ is finitely generated for every $\p \in {\rm End}(S)$.
In view of Corollary~\ref{boup}(ii), we may now obtain~(ii).
\eproof


Similarly, in view of Corollaries~\ref{Cor-ltwo} and \ref{boup}(ii) and Theorem~\ref{fgp}(ii),
we get also

\bt
\label{twobo}
Let $S$ be the monoid (respectively semigroup) defined by a 
finite presentation of the form 
$$\langle A \mid a_1a_2 = a_3a_4\rangle,$$
with $a_1, \dotsc,a_4 \in A$ not necessarily
distinct. Then:
\bi
\item[(i)]
${\rm Per}(\p)$ is finitely generated for every $\p \in {\rm End}(S)$;
\item[(ii)]
for every $\p \in {\rm End}(S)$, there
  exists a constant $R_{\p} > 0$ such that every $a \in {\rm Per}(\p)$
  has period less or equal to~$R_{\p}$. 
\ei
\et

Note that Example \ref{exth} shows that Theorem \ref{twobo} cannot be
generalized to presentations with a relation $u = v$ such that $|u| =
|v| = 3$.

\section*{Acknowledgements}

This work was mostly developed within the activities of FCT's 
(Funda\c c\~ao para a Ci\^encia e a
Tecnologia's) project 
PEst-OE/MAT/UI0143/2013-14 of the 
Centro de \'Algebra da Universidade de Lisboa (CAUL), 
that supported the visit of the first author to CAUL, 
and it was concluded within the FCT project of CEMAT. 

The third author acknowledges support from the European Regional
Development Fund through the programme COMPETE
and the Portuguese Government through FCT under the project 
PEst-C/MAT/UI0144/2013.

\end{document}